\documentclass[10pt]{article}
\usepackage{pb-diagram}
\usepackage{amssymb,latexsym}
\usepackage{amsfonts}
\usepackage{amsmath,amsthm,graphicx}

\usepackage{fontenc,indentfirst, delarray,amsfonts,amsmath,amssymb}
\usepackage{rotating}
\usepackage[T1]{fontenc}

\newcommand{\bd}{\begin{document}}
\newcommand{\ed}{\end{document}}
\newcommand{\bc}{\begin{center}}
\newcommand{\ec}{\end{center}}
\newcommand{\vs}{\vspace}
\newcommand{\hs}{\hspace}
\newcommand{\bq}{\begin{quote}}
\newcommand{\eq}{\end{quote}}
\newcommand{\mb}{\makebox}
\newcommand{\lt}{\left}
\newcommand{\rt}{\right}
\newcommand{\beqa}{\begin{eqnarray*}}
\newcommand{\eeqa}{\end{eqnarray*}}
\newcommand{\beqn}{\begin{eqnarray}}
\newcommand{\eeqn}{\end{eqnarray}}
\newcommand{\bbibl}{\begin{thebibliography}{99}}
\newcommand{\ebibl}{\end{thebibliography}}
\newcommand{\ti}{\times}
\newcommand{\bit}{\begin{itemize}}
\newcommand{\eit}{\end{itemize}}
\newcommand{\ben}{\begin{enumerate}}
\newcommand{\een}{\end{enumerate}}
\newcommand{\lb}{\label}

\newcommand{\hf}{\hspace*{\fill}}
\newcommand{\vf}{\vspace*{\fill}}
\newcommand{\beq}{\begin{equation}}
\newcommand{\eeq}{\end{equation}}
\newcommand{\ba}{\begin{array}}
\newcommand{\ea}{\end{array}}
\newcommand{\del}{\partial}
\newcommand{\ot}{\otimes}
\newcommand{\nn}{\nonumber}
\newcommand{\R}{\mb{$I\!\!R$}}
\newcommand{\C}{{\cal C}}
\newcommand{\M}{{\cal M}}
\newcommand{\E}{{\cal E}}
\newcommand{\N}{{\cal N}}
\newcommand{\B}{{\cal B}}
\newcommand{\X}{{\cal X}}
\newcommand{\Y}{{\cal Y}}
\newcommand{\F}{{\cal F}}
\newcommand{\Rc}{{\cal R}}
\newcommand{\A}{{\cal A}}
\renewcommand{\P}{{\cal P}}
\renewcommand{\S}{{\cal S}}
\newcommand{\es}{\emptyset}
\newcommand{\ci}{\subseteq}
\newcommand{\cs}{\supseteq}
\renewcommand{\u}{\cup}
\renewcommand{\i}{\cap}
\newcommand{\bu}{\bigcup}
\newcommand{\bi}{\bigcap}
\newcommand{\la}{\leftarrow}
\newcommand{\ra}{\rightarrow}
\newcommand{\Ra}{\Rightarrow}
\newcommand{\Lra}{\Leftrightarrow}
\newcommand{\lgra}{\longrightarrow}
\newcommand{\Lgra}{\Longrightarrow}
\newcommand{\lglra}{\longleftrightarrow}
\newcommand{\Lglra}{\Longleftrightarrow}
\renewcommand{\a}{\alpha}
\renewcommand{\b}{\beta}
\newcommand{\g}{\gamma}
\newcommand{\G}{\Gamma}
\renewcommand{\d}{\delta}
\newcommand{\D}{\Delta}
\newcommand{\e}{\varepsilon}
\newcommand{\eps}{\epsilon}
\newcommand{\h}{\eta}
\renewcommand{\l}{\lambda}
\newcommand{\m}{\mu}
\newcommand{\n}{\nu}
\newcommand{\p}{\pi}
\newcommand{\s}{\sigma}
\newcommand{\Si}{\Sigma}
\newcommand{\ta}{\tau}
\newcommand{\ph}{\phi}
\newcommand{\Ph}{\Phi}
\renewcommand{\c}{\chi}
\newcommand{\om}{\omega}
\newcommand{\Om}{\Omega}
\newcommand{\tri}{\triangle}
\newcommand{\rec}[1]{\frac{1}{#1}}
\newcommand{\f}{\frac}
\newcommand{\sm}[2]{\sum_{#1}^{#2}}
\newcommand{\ld}{\ldots}
\newcommand{\ov}{\overline}
\newcommand{\ol}[1]{$\bar{\mb{#1}}$}
\newcommand{\un}{\underline}
\newcommand{\iy}{\infty}

\newcommand{\wt}{\widetilde}
\newcommand{\ds}{\displaystyle}
\newcommand{\bdm}{\begin{displaymath}}
\newcommand{\edm}{\end{displaymath}}

\newcommand{\nin}{\not\in}
\newcommand{\bt}{\begin{tabular}}
\newcommand{\et}{\end{tabular}}
\newcommand{\alter}[2]{\lt\{ \ba {ll}#1 \\ #2 \ea \rt.}
\newcommand{\alt}[4]{\lt\{ \ba{ll}#1 & \mb{if \,\,}#2 \\ #3 & \mb{if
               \,\,}#4 \ea \rt.}
\newcommand{\alto}[6]{ \lt\{ \ba{ll}#1 & \mb{if \,\,}#2 \\ #3 & \mb{if
               \,\,} #4 \\ #5 & \mb{if \,\,}#6 \ea \rt.}
\newcommand{\altero}[5]{\mb{$\lt\{ \ba {ll}#1 & \mb{if \,\,}#2 \\ #3 &
               \mb{if \,\,} #4 \\ #5 & \mb{otherwise} \ea \rt.$}}

\newcounter{cnt1}
\newcounter{cnt2}
\newcounter{cnt3}
\newcommand{\blr}{\begin{list}{$($\roman{cnt1}$)$} {\usecounter{cnt1}
        \setlength{\topsep}{0pt} \setlength{\itemsep}{0pt}}}
\newcommand{\bla}{\begin{list}{$($\alph{cnt2}$)$} {\usecounter{cnt2}
        \setlength{\topsep}{0pt} \setlength{\itemsep}{0pt}}}
\newcommand{\bln}{\begin{list}{$($\arabic{cnt3}$)$} {\usecounter{cnt3}
                \setlength{\topsep}{0pt} \setlength{\itemsep}{0pt}}}
\newcommand{\el}{\end{list}}
\newcommand{\no}{\noindent}
\newtheorem{Thm}{Theorem}[section]
\newtheorem{Lem}[Thm]{Lemma}
\newtheorem{Prop}[Thm]{Proposition}
\newtheorem{Def}[Thm]{Definition}
\newtheorem{Exm}[Thm]{Example}
\newtheorem{Rem}[Thm]{Remark}
\newtheorem{Cor}[Thm]{Corollory}
\renewcommand{\baselinestretch}{1}
\newcommand{\ilim}{\mathop{\varprojlim}\limits}
\newcommand{\dlim}{\mathop{\varinjlim}\limits}

\begin{document}
\title{Equivariant Twisted Cartan Cohomology Theory}
\author{Debasis Sen}
\date{}
\maketitle
\noindent
\begin{abstract}
 In this note we prove an equivariant version of a result of Cartan \cite{car} for equivariant simplicial cohomology with local coefficients.

\end{abstract}
{\bf Keywords: Simplicial sets, group action, local coefficients, Cartan cohomology theory, equivariant twisted cohomology, generalized\\ Eilenberg-MacLane complex.}
\footnote {\bf The author would like to thank CSIR for its support.}

\section{Introduction}

To generalize Sullivan's theory of rational de Rham complexes on simplicial sets to cochain complexes over arbitrary ring of coefficients, Cartan \cite{car} introduced the notion of 'Cohomology theory'. Over the coefficient ring $\mathbb{Z}$, Cartan's result can be described as follows. Recall that a $\emph{simplicial differential}$ $\emph{graded algebra}$ over $\mathbb{Z}$ is a simplicial object in the category DGA of differential graded algebras over $\mathbb{Z}$, so that for each $p\geq 0$ we have a differential graded algebra
$$(A_p^*,\delta)\colon ~ A_p^0\xrightarrow{\delta} A_p^1\xrightarrow{\delta} A_p^2\rightarrow \cdots $$ together with face and degeneracy maps $\partial_i\colon A_{p+1}^*\rightarrow A^*_p$ and $s_i\colon A_p^*\rightarrow A_{p+1}^*$ which are homomorphism of differential graded algebras satisfying the usual simplicial and differential identities. Then a $\emph{cohomology theory}$ in the sense of Cartan is a simplicial differential graded algebra $A$ over $\mathbb{Z}$ such that
\begin{enumerate}
\item each cochain complex $(A_p^*,\delta)$ is exact and $Z^0 A= Ker(A_*^0\xrightarrow{\delta} A_*^1)$ is a simplicially trivial algebra over $\mathbb{Z}$ (here simplicially trivial means that all the face and degeneracy maps are isomorphisms),
\item the homotopy groups $\pi_i(A_*^n)$ of the simplicial set $A_*^n=\{A_p^n\}_{p\geq 0}$ are trivial for all $i,n\geq 0.$
\end{enumerate}

A cohomology theory $A$ determines a contravariant functor from the category of simplicial sets to DGA which assigns to each simplicial set $K$ the differential graded algebra $A(K)=\{Hom(K, A_*^n)\}_{n\geq 0}$, where $Hom(K, A_*^n)$ is the abelian group of simplicial maps $K\rightarrow A_*^n$ and the differential on $A(K)$ is induced from that of $A$. Then Cartan's theorem states that there is a natural isomorphism
$$H^*(A(K))\cong H^*(K;\mathbb{Z}(A)),$$ where $\mathbb{Z}(A)$ is the abelian group $(Z^0A)_0.$

In \cite{hir}, Hirashima generalized Cartan's result for cohomology with local coefficients. Moreover Cartan's theorem was generalized in \cite{mn} for $G$-simplicial sets, i.e for simplicial sets equipped with an action of a group $G$ by simplicial maps. In the equivariant setting, ordinary cohomology of simplicial sets is replaced by Bredon cohomology of $G$-simplicial sets.

Recently, in \cite{ms} the notion of equivariant cohomology with local coefficients, called Bredon-Illman cohomology with local coefficients, has been formulated for $G$-simplicial sets. This is the simplicial version of the equivariant cohomology with local coefficients for $G$-spaces as introduced in \cite{mm}, which generalizes Bredon-Illman cohomology \cite{br}, \cite{ill}. It is therefore reasonable to prove a version of Cartan's theorem on $G$-simplicial sets for equivariant cohomology with local coefficients. In this note we define the notion of equivariant twisted Cartan cohomology theory and prove a version of Cartan's theorem which reduces to the result of \cite{mn} when the local coefficients system is simple, and, to that of \cite{hir} when $G$ is a trivial group. In the present context Cartan's cohomology theory appears as a contravariant functor from the category $O_G$ of canonical orbits to the category of cohomology theory in the sense of Cartan \cite{car}, satisfying some naturality conditions (see Definition \ref{aa} for details). It may be remarked that the equivariant analogue of a point is the orbit of the point and these orbits are precisely the objects of the category $O_G.$

Throughout the paper, when referring to 'model', we will always mean it to be minimal.

This paper is organized as follows. Section 2 is a review of the basic definitions and results that will be used in the sequel. In section 3, we describe the notion of equivariant twisted cohomology and state a classification theorem which is proved in \cite{ms}. In Section 4, we define the notion of equivariant Cartan cohomology theory and prove our main result.

\section{Preliminaries}
In this section we recall some basic definitions and facts about simplicial sets \cite{may} and related topics.
We denote the category of simplicial sets and simplicial maps by $\mathcal{S}$ and the category of
simplicial groups and simplicial group homomorphisms by $\mathcal{SG}.$ Throughout $G$ will denote a discrete group.
\begin{Def}\cite{moore}

Let $B$ be a simplicial set and $\Gamma$  a simplicial group. Then a graded function
 $$\tau\colon B\lgra \Gamma,~~\tau_q\colon B_q\rightarrow \Gamma_{q-1}$$ is called a $\emph{twisting function}$ if it satisfies the following identities:
\begin{equation*}
\begin{split}
 \partial_0(\tau_q(b))&=(\tau_{q-1}(\partial_0b))^{-1}\tau_{q-1}(\partial_1 b),~~ b\in B_q\\
\partial_i(\tau_q(b))&=\tau_{q-1}(\partial_{i+1}b)~~ i>0\\
s_i(\tau_q(b))&=\tau_{q+1}(s_{i+1}b) ~~i\geq 0\\
\tau_{q+1}(s_0 b)&=e_q,~~e_q \mbox{ being the identity of the group }\Gamma_q.
\end{split}
\end{equation*}

\end{Def}

\begin{Def}\lb{tcp}
 Let $B,F$ be simplicial sets, $\Gamma$ a simplicial group which operates on $F$ from the left, and $\tau\colon B\rightarrow \Gamma$ a twisting function.
A $\emph{twisted cartesian}$ $\emph{product (TCP)}$, with fibre $F$, base $B$ and group $\Gamma$ is a simplicial set, denoted by $F\times_{\tau} B$ which satisfies
$$(F\times_{\tau} B)_n=F_n\times B_n$$
and has face and degeneracy operators
\begin{equation*}
\begin{split}
 \partial_i(f,b)=(\partial_i f,\partial_i b),~~ i>0\\
 \partial_0(f,b)= (\tau(b)\partial_0 f, \partial_0 b)~~~\\
 s_i(f,b)=(s_i f, s_i b)~~ ~i\geq 0
\end{split}
\end{equation*}
 \end{Def}

If $B,F$ are Kan complexes then $F\times_{\tau} B$ is also a Kan complex and the canonical projection $p\colon F\times_{\tau} B\rightarrow B$ is a Kan fibration.

For an abelian group $A$ and an integer $n>1$, let $K(A,n)$ denote a minimal Eilenberg MacLane complex of type $(A,n)$. There is a canonical model of $K(A,n)$ for which the $q$-simplices are described as follows. Consider the simplicial abelian group $C(A,n)$ with $q$-simplices $$C(A,n)_q=C^n(\Delta[q]; A),$$ the group of normalized $n$-cochains of the standard simplicial $q$-simplex $\Delta[q]$ \cite{may}.
The face and degeneracy maps of $C(A,n)$ are given as follows. For $\mu\in C(A,n)_q$, $\alpha\in \Delta[q-1]_n$ and $\beta\in \Delta[q+1]_n$
$$\partial_i\mu(\a) =\mu (\delta_i(\a)),~~ s_j\mu(\b) =\mu(\sigma_j(\b)).$$ Here $\delta_i\colon \Delta[q-1]\lgra\Delta[q]$ and $\sigma_j\colon \Delta[q+1]\lgra \Delta[q]$ are the simplicial maps defined by $\delta_i(\Delta_{q-1})=\partial_i\Delta_q,~\sigma_j(\Delta_{q+1})=\Delta_q$, $\Delta_q=(0,1,\cdots,q)$ being the unique non-degenerate $q$-simplex of $\Delta[q]$.

We have a simplicial group homomorphism $$\delta^n\colon  C(A,n)\lgra C(A, n+1)$$ such that $\delta^n c \in C(A, n+1)_q$ is the usual simplicial coboundary of $c\in C(A,n)_q$. Then $$K(A,n)_q = Ker~\delta^n = Z^n(\Delta[q];A)$$ the group of normalized $n$-cocycles. It may be noted that $K(A,n)$ is a minimal one vertex Kan complex.
\begin{Def}
 Let $\pi$ be a group. A $\emph{$\pi$-module}$ is a pair $(A,\phi)$ where $A$ is an abelian group and $\phi\colon \pi\rightarrow Aut(A)$ a group homomorphism. A $\emph{map of $\pi$-modules}$ $f\colon (A,\phi)\rightarrow (A^{\prime},\phi^{\prime})$ is a group homomorphism $f\colon A\rightarrow A^{\prime}$ such that $$f(\phi(x)a)=\phi^{\prime}(x)f(a)$$
for all $x\in \pi$ and $a\in A$. The category of $\pi$-modules is denoted by $\pi$-mod.
\end{Def}
 Let $(A,\phi)\in \pi\mbox{-mod}$. Then $\pi$ acts on the minimal one vertex Kan complex $K(A,n)$ in the following way: $$x\mu=\phi(x)\circ \mu\mbox{ where }\mu\in K(A,n)_q=Z^n(\Delta[q];A),x\in\pi.$$

The notion of a $\emph{generalized Eilenberg Maclane complex}$ appears in \cite{git}, \cite{hir}, \cite{bfg}. Roughly speaking, a generalized Eilenberg Maclane complex is a one vertex minimal Kan complex having exactly two non-vanishing homotopy groups, one of them being the fundamental group. It appears as the total space of a Kan fibration. Gitler \cite{git} used it in the construction of cohomology operations in cohomology with local coefficients. It also plays a crucial role in classifying cohomology with local coefficients \cite{hir}, \cite{bfg}. It may be remarked that a product of Eilenberg Maclane complexes is also sometimes referred to as a generalized Eilenberg Maclane complex.

A generalized Eilenberg Maclane complex can be constructed as follows.
Let $\overline{W}\pi$ denotes the standard $\overline{W}$ construction \cite{may} of a group $\pi$. Let $(A,\phi)$ be a $\pi$-module. We have a twisting function $$\tau(\pi)\colon \overline{W}\pi\rightarrow \pi,\mbox{ where }\tau(\pi)(x_1,\cdots,x_q)=x_1, x_i\in\pi,$$ and $\pi$ is considered as a simplicial group with each component $\pi$ and all the face and the degeneracy maps are identities.
For $n>1$ let
$$L_{\pi}(A,n)= K(A,n)\times_{\tau(\pi)} \overline{W}\pi,$$
where the right hand side is the twisted cartesian product as defined in the Definition \ref{tcp}. Then it is a one vertex minimal Kan complex whose  fundamental group is $\pi$, $n$-th homotopy group is $A$  and all other homotopy groups are trivial. Moreover the action of the fundamental group $\pi$ on the $n$-th homotopy group $A$ is given by $\phi$ \cite{thur}. We have a canonical map $p\colon L_{\pi}(A,n)\rightarrow \overline{W}\pi$, $p(c,x)=x$ for $c\in K, x\in \overline{W}\pi$, which is a Kan fibration.

For a group $G$, $\emph{the category of canonical orbits}$, denoted by $O_G,$ is a category whose  objects are cosets $G/H$, as $H$ runs over subgroups of $G$. A morphism from $G/H$ to $G/ K$ is a $G$-map. Recall that such a morphism determines and is determined by a subconjugacy relation $g^{-1}Hg\subseteq K$ and is given by $\hat{g}(eH)=gK$. We denote this morphism by $\hat{g}$ \cite{br}.

A contravariant functor from $O_G$ to $\mathcal{S}$ (resp.the category of groups or the category of abelian groups) is called an $\emph{$O_G$-simplicial set}$ (resp. $\emph{$O_G$-group}$ or $\emph{abelian $O_G$-group}$). We denote by $O_G\mathcal{S}$, the category of $O_G$-simplicial sets with morphisms being natural transformations of functors.

A morphism $f\colon T\rightarrow S$ of $O_G$-simplicial sets is called an $\emph{$O_G$-Kan fibration}$ if $f(G/H)\colon T(G/H)\rightarrow S(G/H)$ is a Kan fibration for each subgroup $H$ of $G$. Similarly, an $O_G$-simplicial set $T$ is called an $\emph{$O_G$-Kan complex}$ if each $T(G/H)$ is a Kan complex for each subgroup $H\subseteq G$.

We recall the following definition from \cite{mn}.
\begin{Def}
Given an $O_G$-group $\lambda$ and an integer $n\geq 0$, an $O_G$-Kan complex $T$ is called an $\emph{$O_G$-Eilenberg Maclane complex of type}$ $(\lambda,n)$ if each $T(G/H)$ is a $K(\lambda(G/H),n)$ and  $T(\hat{g})\colon T(G/H)\rightarrow T(G/K)$ is the unique simplicial homomorphism induced by the linear map $\lambda(\hat{g})\colon \lambda(G/H)\rightarrow \lambda(G/K)$, $g^{-1}Hg\subseteq K$, such that $T(\hat{g})_n\colon K(\lambda(G/H),n)_n\rightarrow K(\lambda(G/K),n)_n$ is $\lambda(\hat{g})$.
\end{Def}
It is proved in \cite{mn} that any two $O_G$-Eilenberg Maclane complexes of the same type are naturally isomorphic. We denote an $O_G$-Eilenberg Maclane complex of type $(\lambda,n)$ by $K(\lambda,n).$ Let $n>1$ and $\lambda$ be an abelian $O_G$-group. Using the canonical model of an ordinary Eilenberg Maclane complex as described at the beginning of this section, we have a canonical model of $K(\lambda,n)$, given by $K(\lambda,n)(G/H)_q= Z^n(\Delta[q];\lambda(G/H)).$

  Let $\mathcal{C}$ be a category and $T\colon O_G\rightarrow \mathcal{C}$ be a contravariant functor. An $O_G$-group $\underline{\pi}$ $\emph{is said to act on}$ $T$ if we have a group homomorphism $\phi_H\colon \underline{\pi}(G/H)\rightarrow Aut_{\mathcal{C}}(T(G/H))$ for each subgroup $H$ of $G$ such that for any subconjugacy relation $g^{-1}Hg\subseteq K$, $$\phi_H(\underline{\pi}(\hat{g})v)\circ T(\hat{g})= T(\hat{g})\circ\phi_K(v),~v\in \underline{\pi}(G/K).$$ We denote this action simply by $\phi$. Thus we can talk of an action of $\underline{\pi}$ on an $O_G$-simplicial set, an $O_G$-group etc. If $\underline{\pi}$ acts on an abelian $O_G$-group $T$, then we call $T$ a $\emph{$\underline{\pi}$-module}$.

\section{G-simplicial set, Equivariant Twisted \\ Cohomology and its Classification}
In this section we briefly recall the definition of equivariant twisted cohomology and its homotopy classification from \cite{ms}.
Let $G$ be a discrete group. Recall that a $\emph{$G$-simplicial set}$ is a simplicial set $X=\{X_n\}$ such that each $X_n$ is a $G$-set and the face, degeneracy maps commute with this action. A $G$-simplicial set $X$ is called $\emph{$G$-connected}$ if each fixed point simplicial set $X^H$, $H\subseteq G$, is connected.

Let $X$ be a $G$-simplicial set and $\underline{\pi}$ be an $O_G$-group. Let $\Phi X$ denote the $O_G$-simplicial set defined by $\Phi X(G/H)=X^H,~\Phi X(\hat{g})(x)=gx,~x\in X^H,~g^{-1}Hg\subseteq K.$
\begin{Def}
 Let $T$ be an $O_G$-simplicial set and $\Gamma$ a simplicial $O_G$-group. A natural transformation of functors $\tau\colon T\rightarrow \Gamma$ is called an $\emph{$O_G$-twisting function}$ if $\tau(G/H)\colon T(G/H)\rightarrow \Gamma(G/H)$ is an ordinary twisting function for each subgroup $H$ of G.
\end{Def}
\begin{Exm}\lb{Exm}
 Consider the $O_G$-group $\underline{\pi}$ as a simplicial $O_G$-group $\{\underline{\pi}_n\}_{n\geq 0}$ where $\underline{\pi}_n=\underline{\pi}$ for all $n\geq 0$ and face and degeneracy maps are identity natural transformations. Define $$\tau(\underline{\pi})\colon \overline{W}\underline{\pi}\rightarrow \underline{\pi},~\tau(\underline{\pi})(G/H)([x_1,\cdots,x_q])=x_1,$$
where $[x_1,\cdots,x_q]\in \overline{W}\underline{\pi}(G/H)_q,~x_i\in \underline{\pi}(G/H),~1\leq i\leq q.$ It is routine to check that $\tau(\underline{\pi})$ is an $O_G$-twisting function.
\end{Exm}
\begin{Exm}\lb{twist}
Let $X$ be a $G$-connected $G$-simplicial set and $v$ be a $G$-fixed $0$-simplex in $X$. Let $\underline{\pi}X \colon  O_G\lgra \mathcal{G}rp$ be the $O_G$-group defined as follows. For any subgroup $H$ of $G$,
$$\underline{\pi}X(G/H) =  \pi_1(X^H,v)$$ and for a morphism $\hat{g}\colon G/H\lgra G/K,~~g^{-1}Hg\subseteq K$, $\underline{\pi}X(\hat{g})$ is the homomorphism of fundamental groups induced by the simplicial map $g\colon X^K\lgra X^H.$ We regard $\underline{\pi}X$ as an $O_G$-group complex in the trivial way, that is, $\underline{\pi}X(G/H)_n = \underline{\pi}X(G/H)$ for all $n$.
We choose a $0$-simplex $x$ on each $G$-orbit of $X_0$ and a $1$-simplex $\omega_x\in X^{G_x}$ such that $\partial_0\omega_x=x,\partial_1\omega_x=v$. For any other $0$-simplex $y$ on the orbit of $x$ we define $\omega_{y}=g\omega_x$ if $y=gx$. Then it is an easy check that this is well defined and $\omega_{y}\in X^{G_y}_1$. For a $0$-simplex $x\in X^H,$ let $\xi_H(x)=[\overline{\omega}_x]$ be the homotopy class of $\overline{\omega}_x\colon  \Delta[1]\lgra X^H$. Here for any $q$-simplex $\sigma$ of a simplicial set $Y$, $\overline
{\sigma}\colon \Delta[q]\rightarrow Y$ denote the unique simplicial map satisfying $\overline{\sigma}(\Delta_q)=\sigma.$
Define $$\{\kappa (G/H)_n\}\colon X^H\rightarrow \pi_1(X^H,v)$$ by $$\kappa (G/H)_n(y)= \xi_H(\partial_{(0,2,\cdots,n)}y)^{-1}\circ[\overline{\partial_{(2,\cdots,n)}y}]\circ\xi_H(\partial_{(1,\cdots,n)}y)$$
where $y\in (X^H)_n$ and
$$\partial_{(0,2,\cdots,n)}y=\partial_0\partial_2\cdots \partial_n y, \partial_{(2,\cdots,n)}y=\partial_2\cdots \partial_n y, \partial_{(1,2,\cdots,n)}y=\partial_1\partial_2\cdots \partial_n y.$$

 It is standard that $\kappa (G/H)$ is a twisting function on $X^H.$ We verify that
$$\kappa\colon  \Phi X \lgra \underline{\pi}X, ~~G/H \mapsto \kappa (G/H)$$ is  natural.
Suppose $H$ and $K$ are subgroups such that $g^{-1}Hg\subseteq K.$ Let $z\in X^K_n$.
Then $y=gz\in X^H_n.$ Observe that if $x_1,x_2\in X^K_1$ are $1$-simplexes such that $\overline{x_1}\simeq \overline{x_2},$ as  simplicial maps into $X^K$ then $\overline{y_1}\simeq\overline{y_2}$ as simplicial maps into $X^H$ where $y_i =gx_i,~~i=1,2.$ Thus
\begin{equation*}
 \begin{split}
   \kappa (G/H)_n\circ \Phi X(\hat{g})(z)\\
&= \kappa (G/H)_n(y)\\
&= \xi_H(\partial_{(0,2,\cdots,n)}y)^{-1}\circ[\overline{\partial_{(2,\cdots,n)}y}]\circ\xi_H(\partial_{(1,\cdots,n)}y)\\
&= \xi_H(g\partial_{(0,2,\cdots,n)}z)^{-1}\circ[\overline{g\partial_{(2,\cdots,n)}z}]\circ\xi_H(g\partial_{(1,\cdots,n)}z)\\
&= g\xi_K(\partial_{(0,2,\cdots,n)}z)^{-1}\circ g[\overline{\partial_{(2,\cdots,n)}z}]\circ g\xi_K(\partial_{(1,\cdots,n)}z)\\
&= \underline{\pi}X(\hat{g})\circ \kappa (G/K)_n(z).
  \end{split}
\end{equation*}
Thus $\kappa\colon  \Phi X \lgra \underline{\pi}X$ is an $O_G$-twisting function.
\end{Exm}

Let $X$ be a $G$-simplicial set, $\tau \colon  \Phi X\lgra \underline{\pi}$ be an $O_G$-twisting function and $M$ a $\underline{\pi}$-module, given by $\phi$. We define equivariant twisted cohomology of $(X,\tau,\phi)$ as follows.

We denote the category of abelian $O_G$-groups by $\mathcal{C}_G$. We have a cochain complex in the abelian category $\mathcal{C}_G$ defined by
 $$\underline{C}_n(X)\colon O_G\rightarrow Ab,~~G/H\mapsto C_n(X^H;\mathbb{Z}),$$
where
$C_n(X^H;\mathbb{Z})$ is the free
abelian group generated by the non-degenerate $n$-simplexes of $X^H$ and for any morphism $\hat{g}\colon G/H\rightarrow G/K,~~g^{-1}Hg\subseteq K$ in $O_G$, $\underline{C}_n(X)(\hat{g})$ is given by the map
$g_*\colon C_n(X^K;\mathbb{Z})\rightarrow C_n(X^H;\mathbb{Z}),$ induced by the simplicial map $g\colon X^K\lgra X^H.$ The boundary $\partial_n\colon \underline{C}_n(X)\rightarrow \underline{C}_{n-1}(X)$ is a natural
transformation defined by $\partial_n(G/
H)\colon C_n(X^H;\mathbb{Z})\rightarrow C_{n-1}(X^H;\mathbb{Z}),$ where $\partial_n(G/H)$ is the ordinary boundary map of the simplicial set $X^H$. Dualising this chain complex in the abelian
category $\mathcal{C}_G$ we get the cochain complex
$$\{ C^*_G(X;M)=Hom_{\mathcal{C}_G}(\underline{C}_*(X),M),\delta^n\},$$ which defines the ordinary Bredon cohomology of the $G$-simplicial set $X$ with coefficients $M$ \cite{br}. To define the twisted cohomology of the $G$-simplicial set $X$ we modify the coboundary maps as follows
$$\delta^n_{\tau}\colon C_G^{n}(X;M)\rightarrow C_G^{n+1}(X;M),~~f\mapsto \delta^n_{\tau}f$$
where
$$\delta^n_{\tau}f(G/H)\colon C_{n+1}(X^H;\mathbb{Z})\rightarrow M(G/H)$$
is given by
$$\delta^n_{\tau}f(G/H)(x)=(\tau(G/H)_{n+1}(x))^{-1}f(G/H)(\partial_0 x)+\Sigma_{i=1}^{n+1}(-1)^i f(G/H)(\partial_i x)$$
for $x\in X^H_{n+1}$. Note that the first term of the right hand side is obtained by the given action $\phi$. We denote the resulting cochain complex by
$C^*_G(X;\tau, \phi).$
\begin{Def}
The $\emph{$n^{th}$ equivariant twisted cohomology of $(X,\tau,\phi)$}$ is defined as $$H^n_G(X;\tau,\phi)= H_n(C^*_G(X;\tau, \phi)).$$
\end{Def}

 Suppose that $B,F$ are $O_G$-Kan complexes and $\Gamma$ an $O_G$-group complex. Also assume that $B$ is a $\Gamma$-module and $\kappa\colon  B\rightarrow \Gamma$ an $O_G$-twisting function. Then we have the $O_G$-Kan complex $F\times_{\kappa} B$, defined as $$(F\times_{\tau} B)(G/H) = F(G/H)\times_{\tau(G/H)}B(G/H),~(F\times_{\tau} B)(\hat{g})=(F(\hat{g}),B(\hat{g})),$$
 for each object $G/H$ and morphism $\hat{g}\colon G/H\rightarrow G/K$ of the category $O_G.$ We call this $O_G$-Kan complex the $\emph{$O_G$-twisted cartesian product}$ (TCP), with fibre $F$, base $B$, group $\Gamma$ and twisting $\kappa.$ Observe that the second factor projection gives an $O_G$-Kan fibration $p\colon (F\times_{\tau} B)\rightarrow B.$ We view $(F\times_{\tau} B,p)$ as an object in the slice category (cf. \cite{gj}) $O_G\mathcal{S}/B.$

Let $M$ be a $\underline{\pi}$-module with module structure given by $\phi$. For each subgroup $H$ of $G$, define a group homomorphism $$\psi_H\colon \underline{\pi}(G/H)\rightarrow Aut_{\mathcal{S}}(K(M(G/H),n))$$ as follows. For $u\in \underline{\pi}(G/H)$, let $\psi_H(u)$ be the unique simplicial automorphism of $K(M(G/H),n)$ such that $$\phi_H(u)=\psi_H(u)_n\colon K(M(G/H)_n\rightarrow K(M(G/H)_n,~u\in \underline{\pi}(G/H).$$ This defines an action of the $O_G$-group $\underline{\pi}$ on the $O_G$-Kan complex $K(M,n).$ Therefore we can form the $O_G$-Kan fibration $p\colon K(M,n)\times_{\tau(\underline{\pi})} \overline{W}\underline{\pi}\rightarrow \overline{W}\underline{\pi}$, where $\tau(\underline{\pi})$ is the $O_G$-twisting function as described in the Example \ref{Exm}. If we use the canonical model of $K(M,n)$, the total complex of the resulting $O_G$-Kan fibration is denoted by $L_{\phi}(M,n)$. Since any two models of $K(M,n)$ are naturally isomorphic, $K(M,n)\times_{\tau(\underline{\pi})} \overline{W}\underline{\pi}$ is isomorphic to $L_{\phi}(M,n)$ for any model of $K(M,n).$ We call $L_{\phi}(M,n)$ a $\emph{generalized $O_G$-Eilenberg Maclane complex}$. Note that $L_{\phi}(M,n)(G/H)$ is the generalized Eilenberg Maclane complex $$L_{\underline{\pi}(G/H)}(M(G/H),n)= Z^n(\Delta[-];M(G/H))\times_{\tau(\underline{\pi}(G/H))} \overline{W}\underline{\pi}(G/H).$$

The equivariant twisted cohomology $H_G^*(X;\tau,\phi)$ has been classified by the $O_G$-Kan complex $L_{\phi}(M,n)$ in \cite{ms}. This classification result can be described as follows.

Let $X$ be a $G$-simplicial set and $\tau\colon \Phi X \rightarrow \overline{W}\underline{\pi}$ be an $O_G$-twisting function. It determines an $O_G$-simplicial map
 $\theta (\tau)\colon  \Phi X\lgra \overline{W}\underline{\pi}$ defined as,
$$\theta(\tau)(G/H)\colon X_{q}^{H}\lgra \overline{W}\underline{\pi}(G/H)_q,$$
$$x\mapsto [\tau(G/H)_q(x),\tau (G/H)_{q-1}(\partial_{0}x),\cdots,\tau (G/H)_1(\partial_{0}^{q-1}x)].$$
Let $(\Phi X,~L_{\phi}(M,n))_{\overline{W}\underline{\pi}}$ denote the set of liftings of the map $\theta(\tau)$ with respect to $p\colon L_{\phi}(M,n)\rightarrow \overline{W}\underline{\pi},~p(c,g)=g.$
\begin{Def}
 Let $f,~g\in (\Phi X,~L_{\phi}(M,n))_{\overline{W}\underline{\pi}}.$ Then $f$ and $g$ are said to be $\emph{vertically homotopic}$, written $f\sim_v g$, if there is a map $F\colon \Phi X\times \Delta[1]\rightarrow L_{\phi}(M,n)$ of $O_G$-simplicial sets such that for every object $G/H$ of $O_G$, $F(G/H)$ is a homotopy of the simplicial maps $f(G/H),g(G/H)$ and $p\circ F = \theta (\tau)\circ pr_1,$ where $pr_1\colon  \Phi X\times \Delta[1]\lgra \Phi X$ is the projection onto the first factor.
\end{Def}
Observe that $(\Phi X,\theta(\tau))$ and $(L_{\phi}(M,n),p)$ are objects in the slice category $O_G\mathcal{S}/\overline{W}\underline{\pi}$ and $(\Phi X,~L_{\phi}(M,n))_{\overline{W}\underline{\pi}}$ is the set of morphisms in $O_G\mathcal{S}/\overline{W}\underline{\pi}$ from $(\Phi X,\theta(\tau))$ to $(L_{\phi}(M,n),p).$  The category $O_G\mathcal{S}$ is a closed model category \cite{dk} in the sense of Quillen \cite{qui} and hence $O_G\mathcal{S}/\overline{W}\underline{\pi}$
is also a closed model category \cite{gj}, where $(L_{\phi}(M,n),p)$ is a fibrant object and the above notion of vertical homotopy coincides with the abstract homotopy. Therefore $\sim_v$ is an equivalence relation on the set $(\Phi X,~L_{\phi}(M,n))_{\overline{W}\underline{\pi}}$. Let $[\Phi X,~L_{\phi}(M,n)]_{\overline{W}\underline{\pi}}$ denote the set of equivalence classes.
Then the homotopy classification of equivariant twisted cohomology can be stated as follows.
\begin{Thm}\lb{class}
 Suppose $X$ is a $G$-simplicial set and $\tau\colon \Phi X\rightarrow \underline{\pi}$ is an $O_G$-twisting function. Then $$H_G^{n}(X;\tau,\phi)\cong [\Phi X,~L_{\phi}(M,n)]_{\overline{W}\underline{\pi}},~\mbox{for each }n\geq 0.$$
\end{Thm}

\section{Equivariant Twisted Cartan Cohomology \\Theory}


In this final section we formulate an equivariant version of Cartan's Cohomology theory \cite{car} and prove that Bredon-Illman cohomology with local coefficients of a $G$-simplicial set can be computed by the cohomology of a differential graded algebra determined by a given cohomology theory.

We begin with the following equivariant generalization of Cartan Cohomology theory suitable for our purpose.
\begin{Def}\lb{aa}
 An $\emph{equivariant twisted Cartan cohomology theory}$ is a \linebreak sequence $\mathcal{A}=\{A^i\}_{i\geq 0}$ of simplicial abelian $O_G$-groups $A^i$, together with simplicial differentials $\delta^i\colon A^i\rightarrow A^{i+1}$ such that
\begin{enumerate}
\item For each subgroup $H\subseteq G,$ $\mathcal{A}(G/H)=(A^*(G/H)_*,\delta^*(G/H))$ is a simplicial differential graded algebra over $\mathbb{Z}$.
 \item For each $p\geq 0$,
 $$A^0_p\xrightarrow{\delta_p^0} A^1_p\xrightarrow{\delta_p^1} A^2_p\rightarrow\cdots $$
 is an exact sequence in the abelian category $\mathcal{C}_G$ of abelian $O_G$-groups.
 \item The $O_G$-group $\pi_n\circ A^i$ is the zero $O_G$-group, for all $n,i\geq 0.$
 \item  The simplicial abelian $O_G$-group $Z^0\mathcal{A}=ker(A^0\xrightarrow{\delta^0} A^1)$ is simplicially trivial.
\item For each subgroup $H\subseteq G$ and an integer $i\geq 0$ there is a group homomorphism
$$\psi_H^i\colon Aut((Z^0\mathcal{A})_0(G/H))\rightarrow Aut_{\mathcal{SG}}(A^i(G/H))$$ satisfying
\begin{itemize}
 \item $\delta^i\circ \psi_H^i(\alpha)=\psi^{i+1}_H(\alpha)\circ \delta^i,\alpha\in Aut((Z^0\mathcal{A})_0(G/H))~i\geq 0.$
\item If $g^{-1}Hg\subseteq K$, $\alpha\in Aut((Z^0\mathcal{A})_0(G/H)),~\beta\in Aut((Z^0\mathcal{A})_0(G/K))$ such that $\alpha\circ (Z^0\mathcal{A})_0(\hat{g})=(Z^0\mathcal{A})_0(\hat{g})\circ\beta$ then $$\psi_H(\alpha)\circ A^i(\hat{g})=A^i(\hat{g})\circ\psi_K(\beta).$$
\end{itemize}

\end{enumerate}
\end{Def}

\begin{Exm}
For an abelian group $B$ and an integer $n\geq 0$, let $C(B,n)$ denote the simplicial abelian group and $\delta^n\colon C(B,n)\rightarrow C(B,n+1)$ be the simplicial homomorphism as introduced in the Section $1$. Then, for an abelian $O_G$-group $M$, $\mathcal{A}=\{A^i\}_{i\geq 0}$ where $A^n(G/H)= C(M(G/H),n)$ together with the differential $\delta^n$, defines an equivariant twisted Cartan cohomology theory such that $(Z^0\mathcal{A})_0=M.$
\end{Exm}

\begin{Lem}\lb{cont}
 Let $\mathcal{A}\colon A^0\xrightarrow{\delta}A^1\xrightarrow{\delta}\cdots$ be an equivariant twisted Cartan cohomology theory. Then each $A^n$ is contractible as an object of $O_G\mathcal{S}.$
\end{Lem}
\begin{proof}
Consider the abelian $O_G$-simplicial group $Z^n\mathcal{A}$ defined by $Z^n\mathcal{A}(G/H)=Ker(\delta^n(G/H)\colon A^n(G/H)\rightarrow A^{n+1}(G/H))$, $Z^n\mathcal{A}(\hat{g})=A^n(\hat{g})|_{Z^n\mathcal{A}(G/H)}.$
 For an integer $n\geq 0$ and a subgroup $H$ of $G$, we have a short exact sequence
$$0\rightarrow Z^n\mathcal{A}(G/H)\rightarrow A^n(G/H)\rightarrow Z^{n+1}\mathcal{A}(G/H)\rightarrow 0$$
of simplicial abelian groups. Therefore $A^n(G/H)\rightarrow Z^{n+1}\mathcal{A}(G/H)$ is a principal fibration with fibre $Z^n\mathcal{A}(G/H)$ in the category of simplicial sets, and hence a principal twisted cartesian product (PTCP) of type (W) with group complex $Z^n\mathcal{A}(G/H)$ \cite{may}. This PTCP of type (W) is naturally isomorphic to the universal PTCP of type (W), $W(Z^n\mathcal{A}(G/H))\rightarrow \overline{W}(Z^n\mathcal{A}(G/H))$. But $W(Z^n\mathcal{A}(G/H))$ is contractible. The functions
$$h^H_{q-i}\colon W(Z^n\mathcal{A}(G/H))_q\rightarrow W(Z^n\mathcal{A}(G/H))_{q+1},~ 0\leq i \leq q,~q\geq 0,$$
$$h^H_{q-i}(x_q,\cdots ,x_0)=(0^H_{q+1},\cdots,0^H_{i+1},\partial_0^{q-i}x_q\cdots \partial_0 x_{i+1}\cdot x_i,x_{i-1},\cdots x_0),$$
where $x_j\in Z^n\mathcal{A}(G/H)_j,~0\leq j\leq q$ and $0^{H}_{q+1-r}$ is the zero elements of the abelian group $Z^n\mathcal{A}(G/H))_{q+1-r}~0\leq r\leq q-i$, defines a contraction of $W(Z^n\mathcal{A}(G/H))$ which is natural with respect to morphisms of $O_G.$ Hence $A^n(G/H)$ is also contractible and the contraction is natural. Consequently $A^n$ is contractible as object of $O_G\mathcal{S}.$
\end{proof}

Consider an equivariant twisted Cartan cohomology theory $\mathcal{A}=\{A^i\}_{i\geq 0}$. It determines an abelian $O_G$-group $(Z^0\mathcal{A})_0.$ We denote it by $M.$ Given a $G$-simplicial set $X$, an $O_G$-group $\underline{\pi}$, an $O_G$-twisting function $\tau\colon \Phi X\rightarrow \underline{\pi}$, and a $\underline{\pi}$-module structure $\phi$ on $M$, we shall construct a differential graded algebra over $\mathbb{Z}$ whose cohomology will compute the equivariant twisted cohomology of $(X,\phi,\tau)$.

 Note that, by the second condition of the fifth axiom in the Definition \ref{aa}, $A^n$ becomes a $\underline{\pi}$-module by $(\psi\phi)_H = \psi_H \phi_H\colon \underline{\pi}(G/H)\rightarrow Aut_{\mathcal{SG}}(A^n(G/H))$. To see this, observe that for $g^{-1}Hg\subseteq K$, $v\in \underline{\pi}(G/K)$ we have
 $$\phi_H(\underline{\pi}(\hat{g})v)\circ M(\hat{g})=M(\hat{g})\circ\phi_K(v).$$ Therefore taking $\alpha=\phi_H(\underline{\pi}(\hat{g})v), \beta=\phi_K(v)$ in the second condition of the fifth axiom in the Definition \ref{aa}, we get
 $$\psi_H\phi_H(\underline{\pi}(\hat{g})v)\circ A^n(\hat{g})=A^n(\hat{g})\circ\psi_K\phi_K(v).$$

 Consider the $O_G$-twisting function as introduced in the Example \ref{Exm}. We form the $O_G$-Kan fibration $p\colon A^n\times_{\tau(\underline{\pi})} \overline{W}\underline{\pi}\rightarrow \overline{W}\underline{\pi}$ by taking the  $O_G$-twisted cartesian product as described in Section $3$.

The $O_G$-twisting function $\tau\colon \Phi X \rightarrow \overline{W}\underline{\pi}$ determines a map $\theta(\tau)\colon \Phi X \rightarrow \overline{W}\underline{\pi}$ defined by $$\theta(\tau)(G/H)_q(x)=[\tau(G/H)(x),\tau(G/H)(\partial_0x),\cdots,\tau(G/H)(\partial_0^{q-1}x)],~x\in X^H_q.$$
Let $A^n_{\phi}(X;\tau)=\{f\colon \Phi X\rightarrow {A^n\times_{\tau(\underline{\pi})}\overline{W}\underline{\pi}}|~p f=\theta(\tau)\}.$ This set has an abelian group structure by fibrewise addition, fibrewise inversion and the zero section. We define a differential $\overline{\delta}^n\colon A^n_{\phi}(X;\tau)\rightarrow A^{n+1}_{\phi}(X;\tau)$ by
$$(\overline{\delta}^n f)(G/H)(x)=(\delta^n(G/H)c,b),~f\in A^n_{\phi}(X;\tau),x\in X^H,f(x)=(c,b).$$
It is straightforward to check that $\{A^n_{\phi}(X;\tau),\overline{\delta}\}$ is a cochain complex. Furthermore $A^*_{\phi}(X;\tau)$ admits a graded algebra structure induced from the differential graded algebra $\mathcal{A}.$ The zero element of this algebra is given by the trivial lift $\bf 0$, defined by $$\mathbf{0} (G/H)_q(x)=(0_q^H,\theta(\tau)(G/H)_q(x)),$$ where $x\in X^H_q$ and $0_q^H$ is the zero of the abelian group $A(G/H)_q$. As before we use the notation $[\Phi X,~Z^n \mathcal{A}\times_{\tau(\underline{\pi})}\overline{W}\underline{\pi}]_{\overline{W}\underline{\pi}}$ to denote the set of vertical homotopy classes of liftings of $\theta(\tau)$.
\begin{Prop}
With the above notations, we have
 $$H^n(A^*_{\phi}(X;\tau))=[\Phi X,~Z^n \mathcal{A}\times_{\tau(\underline{\pi})}\overline{W}\underline{\pi}]_{\overline{W}\underline{\pi}}.$$
\end{Prop}
\begin{proof}
Clearly $Ker(\overline{\delta}^n)= (\Phi X,~Z^n \mathcal{A}\times_{\tau(\underline{\pi})}\overline{W}\underline{\pi})_{\overline{W}\underline{\pi}}.$
We now show that
$$\mbox{Im}(\overline{\delta}^{n-1})=\{f\in (\Phi X,~Z^n\mathcal{A}\times_{\tau(\underline{\pi})}\overline{W}\underline{\pi})_{\overline{W}\underline{\pi}}|f\sim_v \bf{0}\}.$$ Let $F\colon f\sim_v \bf{0}$. Consider the following left lifting problem in the closed model category $O_G\mathcal{S}/\overline{W}\underline{\pi}$ (\cite{dk}, \cite{gj}).
\[
 \begin{diagram}
 \node{\Phi X}\arrow{e,t}{(0,\theta(\tau)x)}\arrow{s,l}{i_1}\node{A^{n-1}\times_{\tau(\underline{\pi})}\overline{W}\underline{\pi}}\arrow{s,r}{\overline{\delta}^{n-1} }\\
\node{\Phi X \times \Delta[1]}\arrow{e,b}{F}\arrow{ne,t,..}{\tilde{F}}\node{Z^n\mathcal{A}\times_{\tau(\underline{\pi})}\overline{W}\underline{\pi}}
 \end{diagram}
\]
 Here the $O_G$-simplicial set $\Phi X\times \Delta[n]$ is defined by $$(\Phi X\times \Delta[n])(G/H)=X^H\times \Delta[n],~(\Phi X\times \Delta[n])(\hat{g})=(g,id),  n\geq 0.$$ We identify $\Phi X$ with $\Phi X\times \Delta[0]$. The canonical inclusions $\delta_0,\delta_1\colon \Delta[0]\rightarrow \Delta[1]$ (see Section $1$) induce natural inclusions $i_0,i_1\colon \Phi X\rightarrow \Phi X\times \Delta[1]$. Note that $i_1$ is a trivial cofibration and $\overline{\delta}^{n-1}$ is a fibration in $O_G\mathcal{S}/\overline{W}\underline{\pi}$. Hence the above left lifting problem has a solution $\tilde{F}$. Then $\tilde{F} i_0\in A^{n-1}_{\phi}(X;\tau)$ such that $\overline{\delta}^{n-1}(\tilde{F} i_0)=f$. Therefore $f\in \mbox{Im}(\overline{\delta}^{n-1})$.

On the other hand, suppose that $f=\overline{\delta}^{n-1} h$ for $f\in A^n_{\phi}(X;\tau)$ and $h\in A^{n-1}_{\phi}(X;\tau)$. Then clearly $f\in (\Phi X,Z^n \mathcal{A}\times_{\tau(\underline{\pi})}\overline{W}\underline{\pi})_{\overline{W}\underline{\pi}}$. Composing $h$ with first factor projection map, we get a map $h^{\prime}\colon \Phi X\rightarrow A^{n-1}$ of $O_G\mathcal{S}$. But by the Lemma \ref{cont} $A^{n-1}$ is contractible. Let $H\colon \Phi X\times \Delta[1]\rightarrow A^{n-1}$ be a contracting homotopy for the $O_G$-simplicial set $A^{n-1}$. Then define $\tilde{H}\colon \Phi X\times \Delta[1]\rightarrow A^{n-1}_{\phi}(X;\tau)$ by $\tilde{H}(x,t)=(H(x,t),\theta(\tau)x)$. Clearly $\tilde{H}:h\sim_v \bf{0}$ in $O_G\mathcal{S}/\overline{W}\underline{\pi}$. Hence $\overline{\delta}^{n-1}\circ \tilde{H}:f\sim_v \bf{0}$.
This proves the proposition for $n>0.$

For $n=0$, we note that $H^0(A^*_{\phi}(X;\tau))=(\Phi X,~Z^0 \mathcal{A}\times_{\tau(\underline{\pi})}\overline{W}\underline{\pi})_{\overline{W}\underline{\pi}}$ and two elements in the right hand side are homotopic if and only of they are equal.
\end{proof}

 Observe that the fourth axiom of Definition \ref{aa} implies $Z^0\mathcal{A}$ is an $O_G$-Eilenberg Maclane complex of type $(M,0)$ and hence by  induction $Z^n\mathcal{A}$ is an $O_G$-Eilenberg Maclane complex of type $(M,n)$. To justify this, consider the fibration $$A^n(G/H)\rightarrow Z^{n+1}\mathcal{A}(G/H)$$ with fiber $Z^n\mathcal{A}(G/H),~ H\leq G$. As noted in the Lemma \ref{cont}, this is a PTCP with fibre $Z^n\mathcal{A}(G/H)$. Therefore if $Z^n\mathcal{A}(G/H)$ is minimal then so is  $Z^{n+1}\mathcal{A}(G/H).$ But $Z^0\mathcal{A}(G/H),$ being simplicially trivial, is minimal. Hence by induction it follows that $Z^n\mathcal{A}(G/H)$ is minimal for all $n.$ Now applying the homotopy long exact sequence to the above fibration, and using the third axiom of the Definition \ref{aa} together with induction on $n$, we see that $Z^n\mathcal{A}$ is an $O_G$-Eilenberg Maclane complex of type $(M,n)$. Hence it is isomorphic to the canonical model of $K(M,n).$
Therefore $(Z^n\mathcal{A}\times_{\tau(\underline{\pi})}\overline{W}\underline{\pi},p)$ is isomorphic to $(L_{\phi}(M,n),p)$ as objects in the slice category $O_G\mathcal{S}/\overline{W}\underline{\pi}$. So we have,
\begin{equation*}
 \begin{split}
H^n(A_{\phi}(X;\tau))&= [\Phi X, Z^n\mathcal{A}\times_{\tau(\underline{\pi})}\overline{W}\underline{\pi}]_{\overline{W}\underline{\pi}}\\
& \cong [\Phi X,~L_{\phi}(M,n)]_{\overline{W}\underline{\pi}}.
\end{split}
\end{equation*}
It follows from the Theorem $\ref{class}$ that $$H^n(A_{\phi}(X;\tau))\cong H^n_{G}(X;\phi,\tau).$$
Thus we have proved the following theorem.

\begin{Thm}\lb{th}
 Suppose $\mathcal{A}$ is an equivariant twisted Cartan cohomology theory. Then for every $G$-simplicial set $X$ together with an $O_G$-group $\underline{\pi}$, an $O_G$-twisting function $\tau\colon \Phi X\rightarrow \underline{\pi}$ and an action $\phi$ of $\underline{\pi}$ on the abelian $O_G$-group $(Z^0\mathcal{A})_0$ there is an isomorphism of graded algebras $$H_G^*(X;\tau,\phi)\cong H^*(\mathcal{A}(X;\tau,\phi)),$$  where $\mathcal{A}(X;\tau,\phi)$ denote the graded algebra $(A^*_{\phi}(X;\tau),\overline{\delta}).$
\end{Thm}

 It has been shown in \cite{ms} that for a $G$-connected $G$-simplicial set $X$ with a $G$-fixed $0$-simplex, the simplicial version of Bredon Illman cohomology with local coefficients can be interpreted as an equivariant twisted cohomology for $\underline{\pi} = \underline{\pi}X$ and the $O_G$-twisting function $\kappa$ as described in the Example \ref {twist} (cf. Theorem $4.7$ of \cite{ms}). Combining it with the Theorem \ref{th} we have the following result.
 \begin{Thm}
  Suppose $\mathcal{A}$ is an equivariant twisted Cartan cohomology theory. Given any $G$-connected $G$-simplicial set $X$ with a $G$-fixed $0$-simplex and an action $\phi$ of $\underline{\pi}X$ on $(Z^0\mathcal{A})_0$, let $\mathcal{L}$ be the equivariant local coefficients (cf.Definition 3.1 and page 1020, section 3 $\cite{ms}$) determined by the $\underline{\pi}X$-module $(Z^0\mathcal{A})_0$ on $X$. Then
$$ H^*_G(X;\mathcal{L})\cong H^*(\mathcal{A}(X;\kappa,\phi)).$$
 \end{Thm}

\section{Acknowledgement}
The author had many fruitful discussions with Prof. Goutam Mukherjee and would like to thank him for his useful comments.


{\bf Debasis Sen}\\
Indian Statistical Institute, Kolkata-700108, India.\\
e-mail: dsen\_r@isical.ac.in, sen.deba@gmail.com

\end{document}